\documentclass[11pt,reqno]{amsart}

\title[Nonlinear Instability for the Critical QG]{Nonlinear Instability for the Critically Dissipative Quasi-Geostrophic Equation}
\author{Susan Friedlander}
\address{Department of Mathematics,
University of Southern California, 3620 South Vermont Ave., KAP 108
Los Angeles, CA 90089} \email{\tt susanfri@usc.edu}

\author{Nata\v{s}a Pavlovi\'c}
\address{Department of Mathematics, University of Texas at Austin,
1 University Station, C1200, Austin, TX 78712} \email{\tt
natasa@math.utexas.edu}

\author{Vlad Vicol}
\address{Department of Mathematics,
University of Southern California, 3620 South Vermont Ave., KAP 108
Los Angeles, CA 90089} \email{\tt vicol@usc.edu}

\usepackage{amsmath, amsthm, amssymb}
\usepackage{color}
\theoremstyle{plain}
\newtheorem{theorem}{Theorem}[section]
\newtheorem*{definition}{Definition}
\newtheorem{lemma}[theorem]{Lemma}
\newtheorem{proposition}[theorem]{Proposition}

\theoremstyle{definition}
\newtheorem{remark}[theorem]{Remark}

\def\tilde{\widetilde}
\numberwithin{equation}{section}

\def\R{{\mathbb R}}
\def\Z2{{\mathbb Z}^2}
\def\R2{{\mathbb R}^2}
\def\T2{{\mathbb T}^2}

\def\qq{{ q}}
\def\UU{{ U}}
\def\kk{{ k}}
\def\RR{{ R}}

\def\ob{{\omega_B}}

\def\ddiv{\mathop{\rm div} \nolimits}
\def\diam{\mathop{\rm diam} \nolimits}

\newcommand{\norm}[2]{\Vert #1 \Vert_{#2}}

\begin{document}

%%%%%%%%%%%%%%%%%%%%%%%%% THE ABSTRACT %%%%%%%%%%%%%%%%%%%%%%%%%%%%%%%%%%%

\begin{abstract}
We prove that linear instability implies non-linear instability in the energy norm for the critically dissipative quasi-geostrophic equation.
\end{abstract}

%%%%%%%%%%%%%%%%%%%%%%%% Classification and Keywords %%%%%%%%%%%%%%%%%%%%



\maketitle
%%%%%%%%%%%%%%%%%%%%%%%%%% The Main Part %%%%%%%%%%%%%%%%%%%%%%%%%%%%%%%%%%

\section{Introduction}\label{sec:intro}\setcounter{equation}{0}

A fundamental equation in oceanography and meteorology is the 3 dimensional Navier-Stokes equation in the context of a rapidly rotating, density stratified, viscous, incompressible fluid. Both the forces of rotation and stratification impose a tendency toward 2 dimensionality on the 3 dimensional fluid motion, and this leads to approximate and simpler mathematical models. Important non-dimensional parameters are the Ekman number (the strength of the viscous term relative to rotation) and the Rossby number (the strength of the nonlinearity relative to rotation). In many geophysical problems these parameters are very small. A set of approximations based on asymptotic expansions in powers of these small parameters yields an approximate equation for the 3 dimensional pressure known as the general quasi-geostrophic equation with appropriate boundary conditions. Further simplifying assumptions reduce the problem to the study of a 2 dimensional equation which describes the evolution of the temperature field on a surface that bounds the fluid. In the geophysical fluids literature this equation is known as the surface quasi-geostrophic equation. A derivation of this equation and a discussion of its physical relevance can be found, for example, in Pedlosky \cite{Pe}, Salmon \cite{S}, Held at al \cite{HPGS}. The effects of viscosity are incorporated via a boundary layer analysis and a mechanism known as Ekman layer pumping produces the dissipative term in the 2 dimensional quasi-geostrophic equation.

In the mathematical literature this 2 dimensional equation is often called the dissipative quasi-geostrophic equations (QG equation) with the word {\it surface} being omitted since the equation is 2 dimensional. This equation, for an unknown active scalar $\Theta(x,t)$ representing the temperature on the boundary surface, is given by
\begin{align} \label{eq:forcedQG}
\partial_t \Theta+ \UU \cdot \nabla\Theta +(-\Delta)^{\beta}\Theta=f,
\end{align}
where $\UU(x,t)$ is the velocity vector and $f(x)$ is a given external force. The velocity is coupled with the temperature via a stream function $\Psi(x,t)$:
\begin{align}\label{eq:streamtheta}
\Theta = (-\Delta)^{1/2} \Psi = \Lambda \Psi,
\end{align}
and
\begin{align} \label{eq:streamu}
\UU = \nabla^\perp \Psi = (\partial_{x_2} \Psi, -\partial_{x_1} \Psi) = (R_2 \Theta, -R_1 \Theta),
\end{align}
where $R_i$ is the $i^{th}$ Riesz transform. Our analysis of \eqref{eq:forcedQG} - \eqref{eq:streamu}
considers $x$ in the 2 dimensional torus $[0, 2\pi] = \T2$ and $t \in [0, \infty)$.

Both the non-dissipative and the dissipative QG equations have received much attention following seminal article of Constantin et al \cite{CMT}. They observed a number of similar features between the full 3 dimensional Euler and Navier-Stokes equations and the much simpler QG equations in terms of possible formation of singularities. Recent results concerning the dissipative QG equations include \cite{CC,CCW,CV,CW,DD,DP,J1,J2,KN,KNV,M,W} and references therein.

The appropriate power $\beta$ of the Laplacian in the derivation from the general 3D viscous
quasi geostrophic models and Ekman boundary layer analysis is $\beta = 1/2$. Dimensionally
the 2D QG equation with $\beta = 1/2$ is the analogue of the 3D Navier-Stokes equation.
$\beta = 1/2$ is called the critical case. The first results concerning regularity
of solutions to the dissipative QG equation were given in the simpler (but non-physical) subcritical case where $\beta > 1/2$: see, for example, Constantin and Wu \cite{CW}. In the critical case, $\beta = 1/2$, Constantin, Cordoba and Wu \cite{CCW} proved existence of a unique global solution evolving from any initial data that are small in $L^{\infty}$. Very recently, the smallness assumption was removed independently in breakthrough works of Caffarelli and Vasseur \cite{CV} and Kiselev, Nazarov and Volberg \cite{KNV}. In particular, Caffarelli and Vasseur \cite{CV} used harmonic extension to establish regularity of the Leray-Hopf weak solution. On the other hand, Kiselev et al \cite{KNV} proved the global well posedness of the critical dissipative QG equations with periodic $C^{\infty}$ data. Their argument is based on a certain non-local maximum principle for a suitably chosen modulus of continuity.

In the present article we consider the question of nonlinear instability of a steady solution of the forced critical QG equations. We note that the above mentioned references concern the case $f=0$, but in order to ensure the existence of a large class of steady states we must consider the nontrivially forced problem. In particular, we need to reprove certain results that are known to hold for the unforced equations but not in the forced context, namely the nonlocal maximum principle of Kiselev et al \cite{KNV}.

The main result of this paper is that linear instability implies nonlinear Lyapunov instability for $\Theta$, and hence $U$, in the function space $L^2$. Such results connecting linear and nonlinear instability have been proven under certain restrictions for the 2D Euler equations, see Bardos et al \cite{BGS}, Friedlander and Vishik \cite{VF}, and Lin \cite{L}. There the methods utilize a bootstrap technique where closure relies on the special property of conservation of vorticity which is valid for 2D Euler but not for 3D Euler, where the equivalent instability result is still unproven. This property cannot be utilized for the QG equation because the relation between the temperature and the stream function is not equivalent to the relation between the vorticity and the stream function in the 2D Euler equations. In fact this is one reason why it is conjectured that the QG equations might mimic possible singularity development in the 3D fluid equations.

The result that linear instability implies nonlinear instability in $L^2$ for the Navier-Stokes equations in any dimension was proved in Friedlander et al \cite{FPS} (see also the seminal text of Yudovich \cite{Y}). In this case the special ingredient that permits the bootstrap argument to close is the smoothing property of the Laplacian with respect to the nonlinear term. The arguments in \cite{FPS} carry over directly to the subcritical dissipative QG equation (i.e. $\beta > 1/2$) because the dissipative term again smooths the nonlinear term in  \eqref{eq:forcedQG} - \eqref{eq:streamu}. However the case of the critical QG equation is more subtle because the critical dissipative term ($\beta = 1/2$) and the nonlinear term are now of the same order.

Hence to prove linear instability implies nonlinear instability in $L^2$ for the critical dissipative QG equations via the bootstrap argument requires a different special ingredient. The one we use in this article is the existence of a global bound on $\Vert \nabla \Theta(t) \Vert_{L^\infty}$. This result for the unforced critically dissipative QG was proved in \cite{KNV} and a recent preprint of Kiselev and Nazarov \cite{KN} shows that the result also holds for the equation augmented by a dispersion term. The existence of this global bound for the forced equations is proven in Section~\ref{sec:gwp-perturbation}.

We note that the fairly general abstract theorem of Friedlander et al \cite{FSV} may be applied to the critical QG equations - since the spectrum of the linearized operator is discrete (see Section~\ref{sec:lemmas}) and so the spectral gap condition is satisfied - and shows that linear instability implies nonlinear instability in $H^s$, with $s>2$. The novel result of this present paper is to prove instability in the ``physically natural'' energy space $L^2$.

\subsection*{Organization of the paper}

In section~\ref{sec:formulations} we formulate the stability problem in terms of the temperature $\Theta(x,t)$ perturbed about a steady state $\theta_0(x) \in C^{\infty}$. Also in the same section we define nonlinear stability/instability and we state the main instability result, Theorem~\ref{thm:main}. In section~\ref{sec:lemmas} we study the linear operator $L$ for the dissipative QG equations in perturbation form. This operator is elliptic of order $1$, with compact resolvent, and hence its spectrum is purely discrete for $x \in \T2$. We prove certain properties of $L$ that we will use in the bootstrap argument. Then in section~\ref{sec:main} we use this argument to prove Theorem~\ref{thm:main}. In section~\ref{sec:gwp-perturbation} we prove, in the spirit of \cite{KNV}, that the forced equation has a global $C^\infty$ solution and that $\sup_{t\geq 0} \Vert \nabla \Theta(t) \Vert_{L^\infty} <\infty$. This result in used in the bootstrap argument that proves the main theorem.

\subsection*{Acknowledgements}

We thank Hongjie Dong, Alexander Kiselev, Anna Mazzucato, Roman Shvydkoy and Alexis Vasseur for very helpful discussions. The work of S.F. is supported by NSF grant DMS 0803268. The work of N.P. is supported by NSF grant number DMS 0758247 and an Alfred P.~Sloan Research Fellowship.

\section{Notation and formulation of the result}\label{sec:formulations}

Let $\theta_0$ be the temperature of a smooth steady 2D flow with velocity
$\qq_0$, and smooth force $f$, that is we have
\begin{align}
&\qq_0 \cdot \nabla \theta_0 + \Lambda \theta_0 = f \label{eq:steady1}\\
&\qq_0 = (R_2 \theta_0,-R_1\theta_0) \label{eq:steady2}.
\end{align}
Here we consider $\theta_0,\qq_0,f \in C^\infty(\T2)$. We linearize \eqref{eq:forcedQG} about a the steady state $(\theta_0,\qq_0)$ by writing $\Theta(x,t) = \theta_0(x) + \theta(x,t)$ and $\UU(x,t)= \qq_0(x) + \qq(x,t)$. In such a way we obtain an equation that governs the perturbation $\theta$:
\begin{align}
\partial_t \theta = L\theta  + N(\theta)  \label{eq:linearQG},
\end{align}
where the linear operator $L$ is defined by
\begin{align} \label{eq:L}
L\theta = - \qq_0 \cdot \nabla \theta - \qq \cdot \nabla \theta_0 -
\Lambda \theta,
\end{align}
the velocity is coupled with the temperature via
\begin{align}
  \qq = (R_2 \theta, - R_1 \theta) \label{eq:riesz}
\end{align} and
\begin{align} \label{eq:N}
N(\theta) = - \qq \cdot \nabla \theta.
\end{align}For simplicity of the presentation we let $\theta_0,f,\theta$ have zero mean on the torus, and in the following we shall denote $H^s = \{ v \in H^s(\T2): \int_{\T2} v dx = 0\}$, for all $s\geq 0$.
We define a suitable version of stability (the same definition was used, e.g. in \cite{FPS}, \cite{VF}).

\begin{definition}\label{def:stability}
Let $(X,Z)$ be a pair of Banach spaces. A solution $\theta_0$ of
\eqref{eq:steady1}-\eqref{eq:steady2} is called $(X,Z)$ nonlinearly
stable if for any $\rho>0$, there exists $\tilde{\rho}>0$ so that if
$\theta(0) \in X$ and $\Vert{\theta(0)}\Vert_{Z} < \tilde{\rho}$, then we
have
\begin{enumerate}
  \renewcommand{\labelenumii}{\roman{enumii}}
  \item there exists a global in time solution to
  \eqref{eq:linearQG} such that
  $\theta(t) \in C([0,\infty);X)$;
  \item $\Vert{\theta(t)}\Vert_{Z} < \rho$ for
  a.e.~$t\in [0,\infty)$.
\end{enumerate}
An equilibrium $\theta_0$ that is not stable (in the above sense) is
called Lyapunov unstable.
\end{definition}
The Banach space $X$ is the space where a local existence theorem for the nonlinear equations is available, while $Z$ is the space where the spectrum of the linear operator is analyzed, and where the instability is measured. In the case of the critical dissipative QG we let $X$ be the critical Sobolev space $H^1$ (cf.~\cite{CC,CW,DD,J1,J2,M}), while the growth of the perturbation is considered in the energy space $Z=L^2$. Now we are ready to formulate the main result of the present paper.
\begin{theorem}\label{thm:main}
 Suppose that $\theta_0$ is a smooth mean-free steady state solution of
 the critical dissipative QG, i.e., it solves
 \eqref{eq:steady1}-\eqref{eq:steady2}. If the associated linear operator
 $L$, as defined in \eqref{eq:L}, has spectrum in the unstable
 region, then the steady state is $(H^1,L^2)$ Lyapunov nonlinearly
 unstable.
\end{theorem}

\section{Linearized dissipative QG} \label{sec:lemmas}
\setcounter{equation}{0}

The linear operator $L$ defined in \eqref{eq:L} via
\begin{align*}
  L\theta = - \qq_0 \cdot \nabla \theta - \qq \cdot \nabla \theta_0 - \Lambda \theta
\end{align*}is a pseudo-differential operator  with principal symbol \begin{align*}
a(x,\kk) = - |\kk| + i \qq_0(x)
\cdot \kk,\end{align*} which does not vanish on $\T2 \times \Z2\setminus\{0\}$. Therefore $L$ is elliptic of order $1$. Since $q_0, \nabla \theta_0 \in C^\infty$, for large enough $\alpha >0$, we have that $(L-\alpha I)^{-1}$ is a bounded operator from $L^2$ into $H^1$. Moreover, the domain of $L$
\begin{align}
  {\mathcal D}(L) = \{ v\in H^1({\mathbb T}^2), \int_{{\mathbb T}^2} v dx = 0 \}\subset L^2(\T2)
\end{align} is compactly embedded in $L^2$ by Rellich's theorem, so that resolvent $(L-\alpha I)^{-1}$ is a compact operator. Thus $L$ has discrete spectrum.

Let $\mu$ be the eigenvalue of $L$ with maximal positive real part
over $L^2$. Let $\lambda = \rm{Re}\ \mu$ and $\phi \in L^2$ be
the corresponding eigenfunction\footnote{The steady flow $\qq_0 = (\sin m x_2,0)$ gives an example for which the operator $L$ has unstable eigenvalues over $L^2$. This follows from an extension of the analysis in Friedlander and Shvydkoy \cite{FS} to the dissipative equations (see also Meshalkin and Sinai \cite{MS}).}. For a fixed $0<\delta<C_\lambda$, where $C_\lambda>0$ is a constant depending on $\lambda$ to be determined later, we denote by $L_{\delta}$
\begin{align}\label{eq:Ldelta}
L_\delta = L - (\lambda + \delta)I.
\end{align}
The shift ensures that $L_\delta$ generates a bounded $C_0$-semigroup over $L^2$ and that the resolvent set of $L_\delta$ contains the right half plane. The following lemma shows that $L_\delta$ generates an analytic semigroup over $L^2$.

\begin{lemma} \label{lemma:analytic}
Over $L^2$ the operator $L_\delta$ generates an analytic semigroup.
\end{lemma}
The proof of the lemma modifies the proof of \cite[Theorem 7.2.7]{P},
which shows the analyticity of a strongly elliptic operator of order
$2m$ over $L^2$, to the case of the linearized QG operator, which is elliptic of order 1.
\begin{proof}
Define the operator $G$ via
\begin{align}
  G v = \Lambda v + \qq_0 \cdot \nabla v + \RR(v) \cdot \nabla
  \theta_0 + 2 \beta v = - L v + 2\beta v
\end{align}
where we have denoted $\RR(v) = (R_2 v, -R_1 v)$ and $\beta =
\Vert{\nabla \theta_0}\Vert_{L^\infty}$. Since $\qq_0$ is
divergence-free we have that $G$ satisfies G\"arding's
inequality
\begin{align}\label{eq:re}
  {\rm Re}\ (G v,v) \geq \Vert{\Lambda^{1/2}
  v}\Vert_{L^2}^2 + \beta \Vert{v}\Vert_{L^2}^2.
\end{align}
In the above estimate we also used $\Vert{\RR(v)}\Vert_{L^2} \leq
\Vert{v}\Vert_{L^2}$. Similarly, for every $v\in {\mathcal D}(G)$, we have
\begin{align}\label{eq:im}
  |{\rm Im}\ (Gv,v)| \leq |(Gv,v)| \leq \Vert{\Lambda^{1/2}
  v}\Vert_{L^2}^2 + 3\beta \Vert{v}\Vert_{L^2}^2.
\end{align}
Since $v$ is a scalar, it follows from \eqref{eq:re} and
\eqref{eq:im} that the numerical range $S(G)$ (cf.~\cite[pp.~12]{P})
is contained in the set
\begin{align}
  S_{\vartheta_0}=\{ \lambda \in {\mathbb C} : -\vartheta_0 < \arg
  \lambda < \vartheta_0\},
\end{align}
where $\vartheta_0 = \arctan(3) < \pi/2$. Choosing $\vartheta_0 <
\vartheta<\pi/2$ and defining $\Sigma_\vartheta = \{z \in {\mathbb
C} : |\arg z| > \vartheta\}$, we have that there is a constant
$C= C(\vartheta,\vartheta_0)>0$ such that
\begin{align}
  {\rm dist}(z,S(G)) \geq C |z|,\ \mbox{for all}\ z\in
  \Sigma_\vartheta.
\end{align}
We now claim that all real $x<0$ are in the resolvent set $\rho(G)$
of the operator $G$. Recall that $G = -L + 2 \beta I$, and moreover
that the spectrum of the operator $L$ is contained in the half plane
$\{ z\in {\mathbb C}: {\rm Re}\ z \leq \lambda\}$, where $0<\lambda
= {\rm Re}\ \mu$, and $\mu$ is the eigenvalue of $L$ with largest
real part with associated eigenfunction $\phi$. Since $\qq_0$ is
divergence free we also have that
\begin{align}
  \mu \Vert{\phi}\Vert_{L^2}^2 = (L\phi,\phi) = -
  \Vert{\Lambda^{1/2}\phi}\Vert_{L^2}^2 - (\RR(\phi)
  \cdot \nabla \theta_0, \phi),
\end{align}
and by taking real parts this implies that $\lambda \leq
\Vert{\nabla \theta_0}\Vert_{L^\infty} = \beta$; hence the spectrum
of $G$ is contained in the right half plane, proving the claim.

We have hence proven that $\Sigma_\vartheta$ is contained in the the
complement of $\overline{S(G)}$ and has non-empty intersection with
$\rho(G)$; by \cite[Theorem 1.3.9]{P} we have that $\Sigma_\vartheta
\subset \rho(G)$ and for every $z\in \Sigma_\vartheta$ we have the
resolvent estimate
\begin{align}
  \Vert{R(z:G)}\Vert_{\mathcal{L}(L^2)} \leq
  \frac{1}{{\rm dist}(z:\overline{S(G)})} \leq \frac{1}{C |z|}.
\end{align} Therefore $-G$ is the infinitesimal generator of an
analytic semigroup (cf.~\cite[Theorem 2.5.2]{P}) and so $L_\delta =
-G + (2\beta - \lambda - \delta)I$ generates an analytic semigroup
on $L^2$, since it is a bounded perturbation of $-G$.
\end{proof}

Now we state and prove the lemma that will be used in the proof
of our main result, Theorem \ref{thm:main}.

\begin{lemma}\label{lemma:smooth}
For $0\leq \gamma \leq 1$ there exists a constant $C>0$ such
that
\begin{align}\label{eq:smooth}
\Vert{e^{L_\delta t} v}\Vert_{L^2 \rightarrow L^2}
\leq
\frac{C}{t^\gamma}\Vert{v}\Vert_{L^2}^{1-\gamma}
\Vert{\Lambda^{-1}v}\Vert_{L^2}^\gamma,
\end{align}
for all smooth functions $v\in L^2$, where $C = C (\gamma, \delta,\alpha , \theta_0)$.
\end{lemma}
\begin{proof}
 Since $\qq_0$ is divergence free, it is convenient to
  use the operator $A_\alpha$, defined via
  \begin{align}
    A_\alpha v = - \qq_0 \cdot \nabla v - \Lambda v - \alpha v = L_\delta v + \RR(v) \cdot \nabla \theta_0 - (\alpha - \lambda - \delta) v,
  \end{align} where $\alpha>  \max\{\lambda + \delta , C  \Vert \theta_0 \Vert_{H^{2+\epsilon}}^2\}$, $\epsilon>0$, and $C$ is a sufficiently large dimensional constant. We treat $L_\delta$ as a bounded perturbation of $A_\alpha$. The operator $A_\alpha$ is
also elliptic and has discrete spectrum, so by possibly choosing a different $\alpha$,
we have that $A_{\alpha}^{-1} \in {\mathcal L}(L^2)$.

First, we claim that
  \begin{align}\label{eq:claim2}
  \Vert{A_{\alpha}^{-1} \Lambda v}\Vert_{L^2} \leq C \Vert{v}\Vert_{L^2},
  \end{align}
for all smooth $v\in L^2$ with zero mean. In order prove this, denote
$h = A_{\alpha}^{-1} \Lambda v$, which also has zero mean, and observe that \eqref{eq:claim2} is equivalent to
  \begin{align}\label{eq:claim22}
\Vert h \Vert_{L^2} \leq C \Vert \Lambda^{-1} A_\alpha h \Vert_{L^2}.
    \end{align} The definition of $h$ implies that
    \begin{align*}
      (\Lambda^{-1} A_\alpha h, h) = - (\Lambda^{-1}( \qq_0 \cdot \nabla h),h) - \Vert h \Vert_{L^2}^2 - \alpha \Vert  \Lambda^{-1/2} h \Vert_{L^2}^2,
    \end{align*} and therefore
\begin{align}
  \Vert h \Vert_{L^2}^2 + \alpha \Vert  \Lambda^{-1/2} h \Vert_{L^2}^2 \leq \Vert \Lambda^{-1} A_\alpha h \Vert_{L^2} \Vert  h \Vert_{L^2} + | (\qq_0 \cdot \nabla h,\Lambda^{-1} h)|.\label{eq:energy111}
\end{align}
Note that $(\qq_0 \cdot \nabla \Lambda^{-1/2} h,\Lambda^{-1/2} h)=0$ since $\ddiv \qq_0 = 0$. Using Plancherel's theorem, we write this inner product in terms of Fourier coefficients (cf.~\cite{KV} and references therein)
\begin{align}
  (\qq_0 \cdot \nabla h,\Lambda^{-1} h) &= (\qq_0 \cdot \nabla h,\Lambda^{-1} h) - (\qq_0 \cdot \nabla \Lambda^{-1/2} h,\Lambda^{-1/2} h)\notag\\
&= i (2\pi)^2 \sum\limits_{j+k+l=0} \hat{q}_{0j}\cdot k \left(|l|^{-1/2} - |k|^{-1/2}\right) \hat{h}_k |l|^{-1/2} \hat{h}_l.\label{eq:KV1}
\end{align}In the above summation, the Fourier frequencies $j,k,l \in {\mathbb Z}^2 \setminus \{0\}$ because $\qq_0$ and $h$ are mean free, and $\hat{h}_k$ denotes the $k^{th}$ Fourier coefficient of $h$. Since $|l|=|j+k|$ the triangle inequality gives $ | |l|-|k| | \leq |j|$, and therefore
\begin{align*}
  |k| \left| \frac{1}{|l|^{1/2}} - \frac{1}{|k|^{1/2}}\right| \leq |k| \frac{| |l|^{1/2} - |k|^{1/2}|}{|l|^{1/2}|k|^{1/2}}\leq \frac{|j| |k|}{|l|^{1/2}|k|^{1/2} ( |l|^{1/2}+|k|^{1/2})} \leq |j|.
\end{align*}Therefore, by \eqref{eq:KV1} and the Cauchy-Schwartz inequality we have that
\begin{align}
|(\qq_0 \cdot \nabla h,\Lambda^{-1} h)| &\leq C \sum\limits_{j+k+l=0} |j| |\hat{q}_{0j}|  |\hat{h}_k| |l|^{-1/2} |\hat{h}_l|\notag\\
& \leq C \sum\limits_{j\in {\mathbb Z}^2 \setminus \{0\}} |j| |\hat{q}_{0j}| \sum\limits_{l \in {\mathbb Z}^2 \setminus \{0,-j\}} |\hat{h}_{-j-l}| |l|^{-1/2} |\hat{h}_l|\notag\\
&\leq C \Vert h \Vert_{L^2} \Vert \Lambda^{-1/2} h \Vert_{L^2} \sum\limits_{j\in {\mathbb Z}^2 \setminus \{0\}} |j|^{2+\epsilon} |\hat{q}_{0j}| |j|^{-1-\epsilon}\notag\\
&\leq C \Vert h \Vert_{L^2} \Vert \Lambda^{-1/2} h \Vert_{L^2} \Vert \Lambda^{2+\epsilon} \theta_0 \Vert_{L^2}.\label{eq:KV2}
\end{align}We plug this estimate into \eqref{eq:energy111}  and obtain
\begin{align*}
  \frac 12 \Vert h \Vert_{L^2}^2 + (\alpha - C \Vert \Lambda^{2+\epsilon} \theta_0 \Vert_{L^2}^2 ) \Vert  \Lambda^{-1/2} h \Vert_{L^2}^2 \leq \Vert \Lambda^{-1} A_\alpha h \Vert_{L^2}^2
\end{align*}Since $\alpha > C \Vert \theta_0 \Vert_{H^{2+\epsilon}}^2$, the above estimate proves \eqref{eq:claim22} and $A_{\alpha}^{-1} \Lambda \in {\mathcal L}(L^2)$.

Now we prove that for smooth $v \in L^2$ we have
    \begin{align}\label{eq:claim3}
      \Vert{L_{\delta}^{-1} A_\alpha v}\Vert_{L^2} \leq
      C \Vert{v}\Vert_{L^2},
    \end{align}
for a sufficiently large constant $C>0$. The inequality
\eqref{eq:claim3} follows by writing
    \begin{align}\label{temp1}
      L_{\delta}^{-1} A_\alpha v = v + L_{\delta}^{-1}( \RR(v) \cdot
      \nabla \theta_0) - (\alpha - \delta - \lambda) L_{\delta}^{-1} v,
    \end{align}
    and noting that the operator $L_{\delta}^{-1}$ is bounded on
    $L^2$ (cf.~\cite[Lemma 2.6.3]{P}). Together with the boundedness
    of the Riesz-transforms on $L^2$, \eqref{temp1} implies
    \begin{align}
        \Vert{L_{\delta}^{-1} A_\alpha v}\Vert_{L^2} \leq
        \Vert{v}\Vert_{L^2}(1 + C (\Vert{\nabla
        \theta_0}\Vert_{L^\infty} + \alpha -\delta - \lambda)),
    \end{align}
which proves \eqref{eq:claim3} and therefore $L_{\delta}^{-1} A_\alpha \in {\mathcal L}(L^2)$.

In order to conclude the proof of the lemma we use the fact that $L_\delta$ generates an
    analytic semigroup (cf.~Lemma\ref{lemma:analytic}) and therefore
    (cf.~\cite[Theorem 2.6.13]{P}) we have that
    \begin{align} \label{useanalytic}
      \Vert{e^{L_\delta t} v}\Vert_{L^2 \rightarrow L^2}
      = \Vert{L_{\delta}^{\gamma} e^{L_\delta t} L_{\delta}^{-\gamma} v}
      \Vert_{L^2 \rightarrow L^2}
      \leq \frac{C}{t^\gamma}
      \Vert{L_{\delta}^{-\gamma} v}\Vert_{L^2}.
    \end{align}
Now we bound $ \Vert{L_{\delta}^{-\gamma} v}\Vert_{L^2}$
by interpolating (cf.~\cite[Theorem 2.6.10]{P}) as follows
    \begin{align}
      \Vert{L_{\delta}^{-\gamma} v}\Vert_{L^2} = \Vert{L_{\delta}^{1-\gamma} (L_{\delta}^{-1}
      v)}\Vert_{L^2} &\leq C
      \Vert{v}\Vert_{L^2}^{1-\gamma} \Vert{L_{\delta}^{-1}
      v}\Vert_{L^2}^\gamma \notag\\
      & \leq C \Vert{v}\Vert_{L^2}^{1-\gamma}
      \Vert{(L_{\delta}^{-1} A_\alpha) (A_{\alpha}^{-1}
      \Lambda) (\Lambda^{-1} v)}\Vert_{L^2}^\gamma \notag\\
      & \leq C\Vert{v}\Vert_{L^2}^{1-\gamma}
      \Vert{\Lambda^{-1} v}\Vert_{L^2}^\gamma, \label{Ldel-gam}
    \end{align}
where in order to obtain \eqref{Ldel-gam} we used  \eqref{eq:claim3} and \eqref{eq:claim2}.
Now we conclude the proof of the lemma by combining \eqref{useanalytic}  and \eqref{Ldel-gam}.
\end{proof}

\section{Proof of Theorem \ref{thm:main}}
\label{sec:main} \setcounter{equation}{0}

Here we prove Theorem~\ref{thm:main}. In order to do this we must show that the trivial
solution $\theta = 0$ of
\eqref{eq:linearQG} is $(H^1,L^2)$ Lyapunov unstable. With this goal in mind, we consider a family of solutions
$\theta^\varepsilon$ to
\begin{align}
&\partial_t \theta^\varepsilon = L \theta^\varepsilon +
N(\theta^\varepsilon),\label{eq1}\\
  &\theta^\varepsilon|_{t=0} = \varepsilon \phi,\label{eq2}
\end{align} where $\phi$ is as above an eigenfunction of $L$ associated
with the eigenvalue with maximal positive real part $\lambda$. We will prove the following proposition that clearly implies the desired Lyapunov instability result.

\begin{proposition}\label{prop:critical}
There exist positive constants $\bar{C}$ and $\bar{\varepsilon}\leq 1$ such that for every $\varepsilon \in (0,\bar{\varepsilon})$, there exists $T_\varepsilon > 0$ such that
$\norm{\theta^\varepsilon(T_\varepsilon)}{L^2} \geq \bar{C}$.
\end{proposition}

We remark that if $\theta^\varepsilon(x,t)$ solves \eqref{eq1}--\eqref{eq2}, then the function $\Theta^\varepsilon(x,t)= \theta^\varepsilon(x,t) + \theta_0(x)$ solves the forced QG equations \eqref{fqg1}--\eqref{fqg3}, with initial data $\Theta^\varepsilon(x,0) = \theta_0(x) + \varepsilon \phi(x) \in C^\infty(\T2)$. Moreover, in Lemma~\ref{lemma:KNV} of Section~\ref{sec:gwp-perturbation} we prove that the global smooth solution of the forced QG equations satisfies $\Vert \nabla \Theta^\varepsilon(t) \Vert_{L^\infty} \leq C_{0}^{\varepsilon}$ for all $t\geq 0$, where the constant $C_{0}^{\varepsilon}$ depends solely on the $L^\infty$ and $W^{1,\infty}$ norms of the initial data and the force. For $\varepsilon\in (0,1]$, we have $\Vert \Theta^\varepsilon(0) \Vert_{L^\infty} \leq \Vert  \theta_0 \Vert_{L^\infty} + \Vert \phi \Vert_{L^\infty}$, and similarly $\Vert \nabla \Theta^\varepsilon(0) \Vert_{L^\infty} \leq \Vert  \nabla \theta_0 \Vert_{L^\infty} + \Vert \nabla \phi \Vert_{L^\infty}$, which are independent of $\varepsilon$, and therefore there exists a fixed $C_0>0$ such that $\Vert \nabla \Theta^\varepsilon(t) \Vert_{L^\infty} \leq C_{0}$, for all $\varepsilon \in(0,1]$ and for all $t\geq 0$. We refer the reader to the proof of Lemma~\ref{lemma:KNV} for further details. The triangle inequality then implies that by possibly increasing $C_0$ we have
\begin{align}
  \sup_{t\geq 0} \Vert \nabla \theta^\varepsilon(t) \Vert_{L^\infty} \leq C_0 \label{eq:assumption}
\end{align} for all $\varepsilon \in (0,1]$. We will henceforth denote $\theta^\varepsilon$ simply as $\theta$ and will use the analogous notation for $\qq$. All constants in the following are $\varepsilon$-independent.

\begin{proof}[Proof of Proposition~\ref{prop:critical}]
For $R > C_\phi:= \Vert{\phi}\Vert_{L^2}$ to be chosen later, let
$T=T(R,\varepsilon)$ be the maximal time such that
\begin{align}\label{eq:Tdef}
\norm{\theta(t)}{L^2} \leq  \varepsilon R e^{\lambda t},\quad\
\mbox{for}\ t\in[0,T].
\end{align}
Clearly $T\in (0,\infty]$ due to the strong continuity in $L^2$ of
$t\mapsto \theta(t)$ and the chosen initial condition.

Using Duhamel's formula we write the solution of \eqref{eq1}--\eqref{eq2} as
\begin{align}\label{eq:duhamel0}
\theta(t) = e^{Lt}\varepsilon \phi + B(t),
\end{align}where
\begin{align}\label{eq:def:B}
  B(t) = \int_{0}^{t} e^{L(t-s)} N(\theta)(s) \; ds.
\end{align}
First, we shall prove that
\begin{align} \label{Bestimate}
  \Vert{B(t)}\Vert_{L^2} \leq C_1 \left( \varepsilon R e^{\lambda t}\right)^{1 + \gamma/2},
\end{align}where $\gamma \in(0,1)$ and $C_1 = C (C_0, \lambda,\delta,\gamma)>0$ are constants.
To show \eqref{Bestimate}, we rewrite the operator $B$ and then use
Lemma~\ref{lemma:smooth} as follows:
\begin{align}\label{eq:duhamel1}
  \Vert{B(t)}\Vert_{L^2}
  & = \Vert{\int_{0}^{t} e^{(\lambda+\delta)(t-s)}
    e^{L_\delta(t-s)}N(\theta(s)) \; ds }\Vert_{L^2} \nonumber \\
  & \leq \int_{0}^{t} e^{(\lambda+\delta)(t-s)}
    \Vert{e^{L_\delta(t-s)}N(\theta(s))}\Vert_{L^2 \rightarrow L^2} \; ds \notag\\
  & \leq  C \int_{0}^{t} e^{(\lambda+\delta)(t-s)} \frac{1}{(t-s)^{\gamma}}
    \Vert N(\theta(s))\Vert_{L^2}^{1-\gamma}
    \Vert{\Lambda^{-1} N(\theta(s))}\Vert_{L^2}^\gamma \; ds,
\end{align}
where $\gamma\in(0,1)$ is arbitrary, and $C>0$. In order to bound the factor
$\Vert{\Lambda^{-1} N(\theta(s))}\Vert_{L^2}$ we recall the explicit representation (cf.~\cite{CC}) of the nonlinear term
\begin{equation} \label{cc}
  \Lambda^{-1} ( \RR(\theta) \cdot \nabla \theta) = C_n \left( R_1(\theta
  R_2(\theta)) - R_2(\theta R_1(\theta))\right),
\end{equation}
for some dimensional constant $C_n$, and the fact that the Riesz transforms are bounded on $L^2$ and $L^4$, to obtain that
\begin{align}\label{usecc}
  \Vert{\Lambda^{-1} N(\theta(s))}\Vert_{L^2}
\leq C  \Vert{\theta R_i \theta}\Vert_{L^2}
\leq C \Vert{\theta}\Vert_{L^4}^2,
\end{align}
By interpolating, we have
\begin{align}\label{useccc}
  \Vert \theta \Vert_{L^4} \leq C \Vert \theta \Vert_{L^2}^{1/3} \Vert \theta \Vert_{L^8}^{2/3}.
\end{align}
On the other hand by the Gagliardo-Nirenberg inequality and the H\"older inequality we have that
(cf.~\cite{N})
\begin{align}
  \Vert \theta \Vert_{L^8} &= \Vert \theta^{8/3} \Vert_{L^3}^{3/8}  \leq C \Vert \nabla (\theta^{8/3}) \Vert_{L^{6/5}}^{3/8} \notag\\
& \leq C \Vert \theta^{5/3} \nabla \theta \Vert_{L^{6/5}}^{3/8} \leq C \Vert \theta \Vert_{L^2}^{5/8} \Vert \nabla \theta \Vert_{L^\infty}^{3/8}.\label{usecccc}
\end{align}By combining \eqref{usecc} with \eqref{useccc} and \eqref{usecccc} we obtain
\begin{align}
  \Vert{\Lambda^{-1} N(\theta(s))}\Vert_{L^2}^\gamma \leq C \Vert \theta \Vert_{L^2}^{3\gamma/2} \Vert \nabla \theta \Vert_{L^\infty}^{\gamma/2}\label{eq:bound1}
\end{align}
On the other hand, by H\"older's inequality, and the boundedness of the Riesz transforms on $L^2$, we have
\begin{align}\label{eq:bound2}
  \Vert N(\theta) \Vert_{L^2}^{1-\gamma} \leq \Vert \theta \Vert_{L^2}^{1-\gamma} \Vert \nabla \theta \Vert_{L^\infty}^{1-\gamma}.
\end{align}Recall that by \eqref{eq:assumption} we have $\Vert \nabla \theta(t) \Vert_{L^\infty} \leq C_0$, for all $t\geq 0$. Using assumption \eqref{eq:Tdef} and the fact that $0<\delta<C_\lambda = \lambda \gamma/2$, we substitute the bounds \eqref{eq:bound1}
and \eqref{eq:bound2} into \eqref{eq:duhamel1}, to conclude
\begin{align}\label{eq:duhamel2}
\Vert{B(t)}\Vert_{L^2} \leq C_1
  \left(\varepsilon R e^{\lambda t}\right)^{1 + \gamma/2},
\end{align}
for some positive constant $C_1 = C(C_0,\lambda,\delta,q,\gamma)$,
proving \eqref{Bestimate}. The Duhamel formula \eqref{eq:duhamel0} and the bound \eqref{Bestimate}
imply
\begin{align} \label{duh+bound}
  \Vert \theta(t) \Vert_{L^2} \leq C_\phi \varepsilon e^{\lambda t} + C_1
  \left(\varepsilon R e^{\lambda t}\right)^{1 + \gamma/2}.
\end{align}
Observing that $R$ was chosen such that $R > C_\phi$,
it follows that we have the following estimate on the maximal
time $T$:
%defined in \eqref{eq:Tdef}:
\begin{align}\label{eq:Testimate}
\varepsilon e^{\lambda T} \geq \left( \frac{R-C_\phi}{C_1
R^{1 + \gamma/2}}\right)^{2/\gamma}= : C_2>0,
\end{align}
which clearly holds if $T=\infty$. On the other hand, if $T$ is finite, \eqref{eq:Testimate}
is obtained by combining
the continuity of $t \mapsto \norm{\theta(t)}{L^2}$,
\eqref{eq:Tdef} and \eqref{duh+bound} to obtain
\begin{align*}
  \varepsilon R e^{\lambda T} \leq C_\phi \varepsilon e^{\lambda T} +
  C_1 R^{1 + \gamma/2} \varepsilon e^{\lambda T}
\left(\varepsilon e^{\lambda T}\right)^{\gamma/2},
\end{align*}
which, in turn, implies \eqref{eq:Testimate}.
Therefore we have $T \geq T_\varepsilon$, where we defined
\begin{align}\label{eq:Teps}
T_\varepsilon = \frac 1 \lambda \ln \frac{C_2}{\varepsilon}.
\end{align}
To conclude the proof we must find a lower bound on
$\norm{\theta(T_\varepsilon)}{L^2}$. We use Duhamel's formula \eqref{eq:duhamel0}, the
triangle inequality, and \eqref{eq:duhamel2} to obtain
\begin{align} \label{lowbd-Duh}
\norm{\theta(T_\varepsilon)}{L^2} \geq C_\phi \varepsilon e^{\lambda
T_\varepsilon} - C_1 \left(\varepsilon R e^{\lambda
T_\varepsilon}\right)^{1 + \gamma/2}.
\end{align}
Using \eqref{eq:Teps}, with $C_2$ given by \eqref{eq:Testimate}, the lower bound
\eqref{lowbd-Duh} implies
\begin{align*}
\norm{\theta(T_\varepsilon)}{L^2} &\geq C_2 (C_\phi - C_1 R^{1 +  \gamma/2}
\frac{R - C_\phi}{C_1 R^{1 + \gamma/2}})\\
& = C_2 (2 C_\phi - R):=\bar{C}
> 0,
\end{align*}
by choosing $C_\phi < R < 2 C_\phi$. This concludes the proof of the
proposition which, in turn, implies Theorem~\ref{thm:main}.
\end{proof}

\section{Global well-posedness for the forced QG equation} \label{sec:gwp-perturbation}

In this section, by modifying the argument of Kiselev et al \cite{KNV}, we prove that the forced QG equation has a unique global smooth solution. More precisely, we prove the following:
\begin{lemma}\label{lemma:KNV}
Assume that $\Theta_0,f \in C^\infty$ are $\T2$-periodic
functions with zero mean. Then there exists a unique global in time smooth solution of
\begin{align}
  &\partial_t \Theta + U \cdot \nabla \Theta + \Lambda \Theta = f,\label{fqg1}\\
  & U = R(\Theta) = (R_2 \Theta, - R_1 \Theta),\label{fqg2}\\
  & \Theta(0) = \Theta_0.\label{fqg3}
\end{align}Moreover for all $t\geq 0$ we have
\begin{align}\label{eq:lemma:KNV}
  \Vert{\nabla \Theta(t)}\Vert_{L^\infty}
  \leq C_0,
\end{align}
where $C_0 = C_0(\Vert \Theta_0 \Vert_{L^\infty},\Vert \nabla \Theta_0 \Vert_{L^\infty}, \Vert f \Vert_{L^\infty}, \Vert \nabla f \Vert_{L^\infty})$ is a positive constant.
\end{lemma}
The proof of the lemma is in the spirit of \cite{KNV}, but we additionally need to treat the force term, which a-priori could cause growth of the solution. Since $\Theta(t)$ is mean free, it can be shown a-priori that $\Vert \Theta(t) \Vert_{L^p}$, with $2\leq p \leq \infty$, remains bounded for all time. However the same methods do not work for the subcritical quantity $\Vert \nabla \Theta(t) \Vert_{L^\infty}$, and therefore we need to prove the nonlocal maximum principle of \cite{KNV} for the forced QG equation \eqref{fqg1}--\eqref{fqg3}.~This is achieved by suitably choosing a scaling parameter $B$ and making use of the fact that due to periodicity we do not need to consider arbitrarily large length scales.~We note that the scaling parameter $B$ is used only in the modulus of continuity, whereas the solutions to \eqref{fqg1}--\eqref{fqg3} are not space-time rescaled.

\begin{proof}[Proof of Lemma~\ref{lemma:KNV}]
We recall that a continuous, increasing, unbounded, concave function $\omega: [0,\infty) \rightarrow [0,\infty)$, with $\omega(0)=0$ is a modulus of continuity for a function $f$ if
\begin{align}\label{eq:modulus}
  | f(x)  - f(y) | \leq \omega(|x-y|),
\end{align} for all $x,y\in \R2$. The modulus is  \textit{strict} if the strict inequality holds in \eqref{eq:modulus}. We consider a modulus of continuity that also satisfies
$\omega'(0)<\infty$, and $\lim_{\xi \rightarrow 0+} \omega''(\xi) = - \infty$, namely, as in \cite{KNV} we let \begin{align}\omega(\xi) = \begin{cases}
                    \xi - \xi^{3/2}, & 0\leq \xi \leq \delta, \\
                    \delta - \delta^{3/2} + \gamma \log(1 + \frac 14 \log (\xi/\delta) ), & \xi>\delta,
                  \end{cases}\label{eq:ourmodulus}\end{align}
where $\delta>\gamma>0$ are sufficiently small fixed constants.

Since $\Theta_0 \in C^\infty$, there exists a sufficiently large $B>0$ such that $\Theta_0$ has {\it strict} modulus of continuity $\ob(\xi) = \omega(B\xi)$. The scaling parameter $B$ may be chosen as
\begin{align}
  B = C \Vert \nabla \Theta_0 \Vert_{L^\infty} \exp(\exp(C \Vert \Theta_0 \Vert_{L^\infty})),\label{eq:Bcond1}
\end{align}where $C$ is a sufficiently large positive constant. Moreover, since $\omega$ is unbounded, by possibly increasing $B$ we may ensure that
\begin{align}
  A B^2 \geq \Vert \nabla f \Vert_{L^\infty},\label{eq:Bcond2}
\end{align}
where the fixed dimensional constant $A$ is as in \cite[Lemma]{KNV}, and also
\begin{align}
  \frac{\ob(d)}{d} \geq 4\pi \Vert f \Vert_{L^\infty},\label{eq:Bcond3}
\end{align} where $d= \diam(\T2) = 2\pi \sqrt{2}$ will be fixed throughout this section. We fix a $B$ that satisfies \eqref{eq:Bcond1}--\eqref{eq:Bcond3} and recall that the modulus of continuity is given by
\begin{align}
  \omega_B(\xi) = \omega(B\xi)= \begin{cases}
                    B\xi - (B\xi)^{3/2}, & 0\leq \xi \leq \frac{\delta}{B}, \\
                    \delta - \delta^{3/2} + \gamma \log(1 + \frac 14 \log (\frac{B\xi}{\delta}) ), & \xi>\frac{\delta}{B}.
                  \end{cases}\label{eq:ob}
\end{align} Denote $\ob'(\xi) = B \omega'(B\xi)$ and $\ob''(\xi) = B^2 \omega''(\xi)$. We claim that $\ob(\xi)$ is preserved by the evolution \eqref{eq:linearQG}, so that $\Theta$ is a global solution. We extend $\Theta,U,\Theta_0,f$ to $\T2$-periodic functions on $\R2$.

The first step of the proof is to show that if $\Theta(t)$ has strict modulus of continuity $\ob$ for $t\in[0,T]$, then there exists $\tau>0$ such that $\Theta(t)$ has strict modulus of continuity $\ob$ on $t\in[0,T+\tau)$. Since $\Vert  \nabla \Theta(t) \Vert_{L^\infty} < \ob'(0)$, we have that $\Theta(t) \in C^\infty$ for all $t \in [0,T]$, and by the local regularity theorem (cf.~\cite{CC,J1,J2,M}) for some time $\tau>0$ beyond $T$. We must show that by possibly shrinking $\tau$, we have that
\begin{align}\label{eq:extend}
  |\Theta(x,t) - \Theta(y,t) | < \ob(|x-y|)
\end{align}
for all $t\in(T,T+\tau)$ and $x\neq y\in \R2$.

Define the compact set $K = [-2\pi,2\pi]^2 \times [-2\pi,2\pi]^2 \subset {\mathbb R}^4$. Since $\Theta$ is $\T2$-space periodic, we have that for any $(x,y) \in {\mathbb R}^4$, with $x\neq y$, there exist $(x',y') \in K$, with $x'\neq y'$, such that $|x'-y'| \leq |x-y|$, $\Theta(x,t) = \Theta(x',t)$, and $\Theta(y,t) = \Theta(y',t)$. Because $\ob$ is increasing, if \eqref{eq:extend} holds for all $(x',y')\in K$ with $x'\neq y'$, then we have that for all $x\neq y \in \R2$
\begin{align*}
  |\Theta(x,t) - \Theta(y,t)| = |\Theta(x',t) -\Theta(y',t)|  < \ob(|x'-y'|) \leq \ob(|x-y|).
\end{align*}

Therefore it is sufficient to prove that there exists $\tau>0$ such that \eqref{eq:extend} holds for $x\neq y$, with $(x,y) \in K$. By assumption, there exists $\epsilon >0$ such that $\Vert \nabla \Theta(T) \Vert_{L^\infty}< \ob'(0) - 2 \epsilon$, and by continuity, for small enough $\tau$ we have that $\Vert \nabla \Theta(t) \Vert_{L^\infty} < \ob'(0)- \epsilon$ for all $t\in[T,T+\tau)$. Therefore for $(x,y)\in K$, with $0<|x-y|= \xi <\rho$, where $\rho \leq \min(\delta/B,\epsilon^2/B^3)$, we have
\begin{align*}
  | \Theta(x,t)-\Theta(y,t) | \leq \xi \Vert \nabla \Theta(t) \Vert_{L^\infty} < \xi(B-\epsilon) \leq B\xi - (B\xi)^{3/2} = \ob(\xi),
\end{align*} for all $t\in[T,T+\tau)$. On the other hand, due to the continuity in time of $|\Theta(x,t) - \Theta(y,t)|$, the compactness of the set $\{ (x,y)\in K: |x-y|\geq \rho\}$, and and the fact that \eqref{eq:extend} holds at $t=T$, we have that there is a sufficiently small $\tau>0$ such that \eqref{eq:extend} holds for all $(x,y)\in K$, $x\neq y$, and $t\in[T,T+\tau)$.

The second part is to rule out the case in which there exists $T>0$ and $x\neq y \in \R2$ such that $\Theta(x,T) - \Theta(y,T) = \ob(|x-y|)$ (cf.~\cite{KNV}). Note that by the periodicity of $\Theta$, for such $x\neq y\in \R2$ fixed, there exist $x',y' \in \T2$ such that
\begin{align*}
  \ob(|x-y|) &= \Theta(x,T)-\Theta(y,T)\\
  & = \Theta(x',T)-\Theta(y',T) \leq \ob(|x'-y'|) \leq \ob(d),
\end{align*} and since $\ob$ is increasing, we must have $0< \xi = |x-y| \leq d=\diam(\T2)$. We conclude by showing that
\begin{align}\label{eq:obstruction}
  \tfrac{d}{dt} ( \Theta(x,t) - \Theta(y,t))|_{t=T} < 0,
\end{align}contradicting the fact that the {\it strict} modulus of continuity is lost at $t=T$. In the following we suppress the time dependence of $\Theta$ and $U$, since we work at $t=T$ fixed.

Since $\Theta$ has modulus of continuity $\ob(\xi)$, we know (cf.~\cite[Lemma]{KNV}) that $U$ has modulus of continuity $\Omega_B(\xi)$, where we defined
\begin{align*}
  \Omega_B(\xi) = A\left( \int_{0}^{\xi} \frac{\ob(\eta)}{\eta} d\eta + \xi \int_{\xi}^{\infty} \frac{\ob(\eta)}{\eta^2} d\eta\right),
\end{align*} for some positive constant $A$. Then as in \cite[Section 4]{KNV} we have that
\begin{align}
  |(U \cdot \nabla \Theta)(x) - (U \cdot \nabla &\Theta)(y)| \leq \left|\lim_{h\rightarrow 0^+} \frac{\ob(\xi + h |U(x)-U(y)|) - \ob(\xi)}{h}\right| \notag\\
&\leq |U(x)-U(y)| \ob'(\xi) \leq \Omega_{B}(\xi) \ob'(\xi).\label{eq:T1}
\end{align}
The dissipative terms are estimated as in \cite[Section 5]{KNV}, namely by the negative quantity
\begin{align}\label{eq:dissipativebound}
  M_B(\xi) &= \frac{1}{\pi} \int_{0}^{\xi/2} \frac{\ob(\xi + 2\eta) + \ob(\xi-2\eta) - 2\ob(\xi)}{\eta^2} d\eta \notag\\&\qquad + \frac{1}{\pi} \int_{\xi/2}^{\infty} \frac{\ob(2\eta + \xi) - \ob(2\eta-\xi) - 2 \ob(\xi)}{\eta^2} d\eta.
\end{align}
Lastly, the force term is estimated using the mean value theorem
\begin{align}
  |f(x) - f(y)| \leq F_B(\xi) = \begin{cases}
                    \xi\ \Vert \nabla f \Vert_{L^\infty}, & 0\leq \xi \leq \frac{\delta}{B}, \\
                    2\ \Vert f \Vert_{L^\infty}, & \xi>\frac{\delta}{B}.
                  \end{cases}\label{eq:fbound}
\end{align}Thus, in order to conclude the proof of \eqref{eq:obstruction}, we must show that for all $0<\xi\leq d$, we have
\begin{align}
  \Omega_B(\xi) \ob'(\xi) + F_B(\xi) + M_B(\xi) < 0.\label{eq:toprove}
\end{align}
First we treat the case $0<\xi \leq \delta/B$. By keeping track of $B$, and using condition \eqref{eq:Bcond2}, similar arguments as in \cite[Section 7]{KNV} show that
\begin{align*}
  \Omega_B(\xi) \ob'(\xi) + F_B(\xi) + M_B(\xi) &\leq A B^2 \xi(3+ \log\frac{\delta}{B\xi}) + \xi \Vert \nabla f \Vert_{L^\infty} + \frac{\xi}{\pi} \ob''(\xi)\\
  &\leq B^2 \xi \left( A (4 + \log\frac{\delta}{B\xi}) - \frac{3}{4\pi} (B\xi)^{-1/2}\right).
\end{align*}Since we have $0<B\xi\leq \delta$, the above quantity is strictly negative if $\delta$ is sufficiently small. Note that $\delta$ does not depend on $B$.

For the case $\delta/B \leq \xi \leq d$, we follow the estimates in \cite[Section 8]{KNV} to conclude that if $\gamma$ and $\delta$ are sufficiently small, independent of $B$, then
\begin{align*}
  \Omega_B(\xi) \ob'(\xi) + F_B(\xi) + M_B(\xi) &\leq A\gamma \frac{\ob(\xi)}{\xi} + 2\Vert f \Vert_{L^\infty} - \frac{1}{\pi} \frac{\ob(\xi)}{\xi}.
\end{align*}But $B$ was chosen so that \eqref{eq:Bcond3} is satisfied, i.e. $ 2 \Vert f \Vert_{L^\infty} \leq \ob(d)/2\pi d$. Because on $[\delta/B,\infty)$ the function $\ob(\xi)/\xi$ is decreasing, for any $\xi \in [\delta/B,d]$ we have that $2 \Vert f \Vert_{L^\infty} \leq \ob(d)/2\pi d \leq \ob(\xi)/2\pi\xi$. Thus
\begin{align*}
  A\gamma \frac{\ob(\xi)}{\xi} + 2\Vert f \Vert_{L^\infty} - \frac{1}{\pi} \frac{\ob(\xi)}{\xi} \leq \left(A\gamma + \frac{1}{2\pi} - \frac{1}{\pi}\right) \frac{\ob(\xi)}{\xi} < 0,
\end{align*} if $\gamma$ is sufficiently small, independent of $B$. Therefore \eqref{eq:toprove} holds for all $0<\xi\leq d$, and so \eqref{eq:obstruction} is proven. Therefore the solution $\Theta(t)$ exists for all time and has strict modulus of continuity $\ob$, which implies that $\Vert \nabla \Theta(t) \Vert_{L^\infty} < \ob'(0) = B$ for all $t\geq 0$, concluding the proof of the lemma.
\end{proof}

\begin{remark}
  We note that it is also possible to adapt the De Giorgi-type techniques used by Caffarelli and Vasseur \cite{CV} to treat the forced QG equation. First one proves boundedness of the solution in $L^2$ using energy estimates, and then the similarly to \cite[Section 2]{CV} one obtains boundedness (not decay) for all time of $\Theta(t)$ in $L^\infty$ and of $U(t)$ in $BMO$. The second step is to show that the solution is actually H\"older and that it remains bounded in this space for all $t\geq 0$, i.e.~adding a smooth force does not create additional difficulties. Since this is already subcritical regularity, in the third step it is standard to bootstrap to higher regularity and prove that the $W^{1,\infty}$ norm of $\Theta(t)$ is bounded in time.
\end{remark}

\end{document}